\documentclass[11pt]{amsart}
\title{Towards commutator theory for relations. III}
\author{Paolo Lipparini}
\address{Dipartimento di Matematica,
 Viale della Ricerca Scientifica,
II Universitade di Roma (Tor Vergata),
 ROME 
ITALY}
\thanks{The author has received support from MPI and GNSAGA} 

\email{lipparin@axp.mat.uniroma2.it}
\urladdr{http://www.mat.uniroma2.it/\textasciitilde lipparin}
\newtheorem{Theorem}{Theorem}

\newtheorem{thm}[Theorem]{Theorem}

\theoremstyle{definition}

\newcommand{\alg}{\mathbf} 
\def\v{\mathcal V}  

\newcommand{\spacex}{$ $} 

\newcommand{\adma}{\mathrm{Adm({\alg A})}} 
\newcommand{\admb}{\mathrm{Adm({\alg B})}} 
\newcommand{\admx}{\mathrm{Adm({\alg X})}} 

\newcommand{\cona}{\mathrm{Con({\alg A})}} 

\begin{document}

\begin{abstract} 
We derive consequences from the existence of a
term which satisfies Mal'cev identities (characterizing 
permutability) modulo two functions $F$ and $G$ 
from admissible relations to admissible relations.
We also provide characterizations of varieties 
having a Mal'cev term modulo $F$ and $G$.
\end{abstract} 

\maketitle
\bigskip 

Given an algebra $ {\alg A} $, let   $ \adma $ denote
the set of all reflexive and admissible relations on $ {\alg A} $
(we shall use the words ``admissible'' and ``compatible'' interchangeably).

If $ {\alg A} $ is an algebra, and  
$F: \adma \to \adma$, 
$G: \adma \to \adma$, 
we say that a ternary term $t$ of $ {\alg A} $ is
{\em Mal$'$cev modulo $F$ and $G$}
if and only if  
\[
a F(R) t(a,b,b) \qquad \text{ and } \qquad
t(a,a,b) G(R) b,
\]
whenever $a, b \in {\alg A} $, 
$R \in \adma$,
and $aRb$.
An alternative name for the above notion is
a {\em weak difference term modulo $F$ and $G$}, 
or simply an {\em $F$-$G$-difference term}:
we used this terminology in \cite[p. 199]{lptams}, in the case when
$F: \cona \to \cona$, 
$G: \cona \to \cona$.

For a relation $R$ on some algebra, let $R^\circ$ denote the smallest
tolerance containing $R$, and let $R^-$ denote the converse
of $R$. $R^*$ is the transitive closure of $R$, and
$Cg(R)$ is the smallest congruence containing $R$.
$\overline{R}$ denotes the least compatible relation
containing $R$. 
Notice that if $R$ and  $S$ are compatible, then 
$R \circ S$ is compatible, too.
Intersection is denoted by juxtaposition.
$R \circ_n S$ is $R \circ S \circ R \circ S \dots $
with $n-1$ occurrences of $\circ$.
$R^n= R \circ_n R$. By convention, we put
$R^0= 0$, where $0$ denotes the identity relation
(the smallest reflexive relation).
$R+S= \bigcup_{n} R \circ_n S$.

\begin{thm} \label{wdta} 
Suppose that
$ {\alg A} $ has a term 
Mal$'$cev modulo $F$ and $G$.
Then for all
reflexive admissible relations $R,S, R_1, R_2, \dots \in \adma$,
and for  arbitrary relations $ \theta, \theta_1, \theta_2 \subseteq {\alg A}^2$,
the following hold:

(i) If $a R b \theta_1 c$,  $a \theta_2 d S c$, and $b \theta d$ (that is,  the situation 
 pictured in the following diagram occurs):
\[
\begin{array}{rcl} 
                              &b&  \\
  ^R \nearrow      &\big|  &   \searrow ^{\theta_1}\\
a          \ \quad  & \ \,\big| {\theta} & \quad \ c\\
  _{\theta_2} \searrow     &\downarrow &  \nearrow _S\\
                                & d &\\
\end{array}
\] 
then 
$(a,c) \in  F(R) \circ \overline{\theta_2 \cup \theta \cup \theta_1 } \circ G(S ) $.

(ii) $R \circ \theta \circ S \subseteq 
F(R) \circ \overline{(R \circ \theta)\cup (\theta \circ S)} \circ G(S)$. 

(iii) $R  \circ S \subseteq 
F(R) \circ \overline{R \cup S} \circ G(S)
\subseteq F(R) \circ S \circ R \circ G(S)$. 

(iv)
$R \circ_{n+2} S \subseteq 
F(R)\circ \left(\overline{F(R) \cup F(S)}\right)^n  \circ \overline{R \cup S} 
\circ \left(\overline{G(R) \cup G(S)}\right)^n \circ G(S^\bullet)$, 
for every $n\geq 0$, 
where $ S^\bullet =S$ if $n$ is even,
and $ S^\bullet =R$ if $n$ is odd.

(v)
$R + S \subseteq 
(F(R)+ F(S))  \circ \overline{R \cup S} \circ (G(R)+ G(S))\subseteq
(F(R)+ F(S))  \circ R \circ S \circ (G(R)+ G(S))
$.

(vi) $R \circ R \subseteq F(R) \circ R \circ G(R)$. 

(vii) $R ^{n+1} \subseteq (F(R))^n \circ R \circ (G(R))^n$, 
for every $n\geq 0$.

(viii) $R ^* \subseteq (F(R))^* \circ R \circ (G(R))^*$.

(ix) $R^- \subseteq F(R^-) \circ R \circ G(R^-)$; and
$R \subseteq F(R) \circ R^- \circ G(R)$.

(x) $R + S^- \subseteq 
(F(R)+ F(S^-))  \circ \overline{R \cup S} \circ (G(R)+ G(S^-))\subseteq
(F(R)+ F(S^-))  \circ R \circ S \circ (G(R)+ G(S^-))
$.

(xi) $Cg(R) \subseteq  (F(R)+ F(R^-))  \circ R  \circ (G(R)+ G(R^-))$. 

(xii) For $n\geq 2$, $R_1  \circ R_2 \circ \dots \circ R_n \subseteq$
\begin{center}
$F(R_1) \circ \overline{F(R_1) \cup F(R_2)} \circ 
\overline{F(R_1) \cup F(R_2) \cup F(R_3)} \circ
 $\\
$\dots \circ \overline{F(R_1) \cup F(R_2) \cup \dots 
\cup F(R _{n-2} )} \circ
\overline{F(R_1) \cup F(R_2) \cup \dots 
\cup F(R _{n-1} )} \circ
$\\$
\overline{R_1 \cup R_2 \cup \dots \cup R _{n-1}  \cup R _{n}} \circ$\\
$\overline{G(R_2) \cup G(R_3) \cup \dots 
\cup G(R _{n} )} \circ
\overline{G(R_3) \cup G(R_4) \cup \dots 
 \cup G(R _{n} )} \circ$\\
$\dots \circ
\overline{ G(R _{n-2} ) \cup G(R _{n-1} ) \cup G(R _{n} )} \circ
\overline{G(R _{n-1} ) \cup G(R _{n} )} \circ
G(R_n)$
\end{center} 

(xiii) $R_1  + R_2 + \dots + R_n \subseteq 
\big(
F(R_1) + F(R_2) + \dots + F(R _{n-1} ) + F(R _{n} )
\big) \circ
\overline{R_1 \cup R_2 \cup \dots \cup R _{n-1}  \cup R _{n}} \circ
\big(
G(R_1) + G(R_2) + \dots + G(R _{n-1} ) + G(R _{n} )
\big)$. 

(xiv) $Cg(R_1 \cup R_2 \cup \dots \cup R_n) \subseteq$ 
\begin{center}
$\big(
F(R_1) + F(R _{1}^-) + F(R_2) + F(R _{2}^-) +\dots 
+F(R _{n} ) + F(R _{n}^-)
\big) \circ
$\\$
\overline{R_1 \cup R_2 \cup \dots \cup R _{n-1}  \cup R _{n}} \circ
$\\$
 \big(
G(R_1) + G(R _{1}^-) + G(R_2) + G(R _{2}^-) +\dots 
+G(R _{n} ) + G(R _{n}^-)
\big)
$
\end{center}
\end{thm} 

\begin{proof}
(i) 
$a  F(R) t(a,b,b) \overline{\theta_2 \cup \theta \cup \theta_1 } t(d,d,c) G(S ) c$.

 (ii) If $a R b \theta d S c$ then
$a R b \theta d $ and $ b \theta d S c$. Letting
$ \theta_2= R \circ \theta  $ and 
$ \theta_1= \theta  \circ S$ we get the conclusion from (i),
noticing that $\theta \subseteq R \circ \theta$, since $R$ is reflexive.

The first inclusion in (iii) is the particular case $\theta=0$ of (ii).
The second inclusion is trivial, since $R$ and $S$ are reflexive, and
$S \circ R$ is compatible. 

(iv) is proved by induction on $n$.
The base $n=0$ is given by (iii).

Suppose that (iv) holds for some $n$, and that
$(a,c) \in R \circ_{n+3} S $. 
This means that there are $b$, $d$ such that 
$a R b (S \circ_{n+1} R) d R^\bullet c $,
where $ R^\bullet =R$ if $n$ is even,
and $ R^\bullet =S$ if $n$ is odd.
This implies $a (R  \circ_{n+2} S) d $, and
$ b (S \circ_{n+2} R) c $.
Letting 
\[ 
 \theta_2= F(R)\circ \left(\overline{F(R) \cup F(S)}\right)^n  \circ \overline{R \cup S} 
\circ \left(\overline{G(R) \cup G(S)}\right)^n \circ G(S^\bullet)
\]
we get, by the inductive assumption, 
$R \circ_{n+2} S \subseteq \theta_2$, hence
$a \theta_2 d$.
Symmetrically, letting
\[ 
 \theta_1= F(S)\circ \left(\overline{F(R) \cup F(S)}\right)^n  \circ \overline{R \cup S} 
\circ \left(\overline{G(R) \cup G(S)}\right)^n \circ G(R^\bullet)
\] 
we get $b \theta_1 c$
(notice that $\cup$ is a commutative operation).

Letting $ \theta= S \circ_{n+1} R $, we have 
$ b \theta d $. Notice that the inductive assumption also gives
$ \theta= S \circ_{n+1} R  \subseteq 
R \circ_{n+2} S \subseteq \theta_2$.
Thus we get from clause (i)
$(a,c) \in  F(R) \circ \overline{\theta_2 \cup \theta \cup \theta_1 } \circ G(R^\bullet ) $,
which proves the result for $n+1$, since
\begin{multline*}
\theta_2 \cup \theta \cup \theta_1 \subseteq
\\
(F(R)\cup F(S))\circ \left(\overline{F(R) \cup F(S)}\right)^n  \circ \overline{R \cup S} 
\circ \left(\overline{G(R) \cup G(S)}\right)^n \circ (G(S^\bullet) \cup G(R^\bullet))
 \end{multline*} 
hence
\begin{multline*}
\overline{ \theta_2 \cup \theta \cup \theta_1} \subseteq
\\
\overline{F(R)\cup F(S)} \circ \left(\overline{F(R) \cup F(S)}\right)^n  \circ \overline{R \cup S} 
\circ \left(\overline{G(R) \cup G(S)}\right)^n \circ 
\overline{G(R) \cup G(S)}=
\\
\left(\overline{F(R) \cup F(S)}\right)^{n+1}  \circ \overline{R \cup S} 
\circ \left(\overline{G(R) \cup G(S)}\right)^{n+1}
\end{multline*} 

Indeed, from the above identities, we get
$$
R \circ _{n+3} S \subseteq 
F(R) \circ
\left(\overline{F(R) \cup F(S)}\right)^{n+1}  \circ \overline{R \cup S} 
\circ \left(\overline{G(R) \cup G(S)}\right)^{n+1}
\circ G(R^\bullet)
$$
which completes the induction step.

(v) is immediate from (iv).

(vi), (vii) and (viii) are the particular cases $S=R$ of, respectively, 
 (iii), (iv) and (v).

(ix) If $a R^- b$, that is, $b R a$, then 
$a  F(R^-) t(a,b,b) R t(a,a,b) G(R^-) b$.
The second formula follows from the first one, applied with $R^-$
in place of $R$, since $R^{--}=R$.
 
(x)
By (ix) with $S$ in place of $R$   
we get 
$S^- \subseteq F(S^-) \circ S \circ G(S^-)$.
Hence, 
$R \cup S^- \subseteq F(S^-) \circ (R \cup S) \circ G(S^-)$, and
$ \overline{ R \cup S^-} \subseteq F(S^-) \circ \overline{R \cup S} \circ G(S^-)$,
since $F(S^-)$ and $G(S^-)$ are compatible.

By (v) with $S^-$ in place of $S$, we get   
$R + S^- \subseteq 
(F(R)+ F(S^-))  \circ \overline{R \cup S^-} \circ (G(R)+ G(S^-))
\subseteq 
(F(R)+ F(S^-))  \circ 
F(S^-) \circ \overline{R \cup S} \circ G(S^-)
 \circ (G(R)+ G(S^-))
=
(F(R)+ F(S^-))  \circ 
\overline{R \cup S} 
 \circ (G(R)+ G(S^-))
$.

(xi) is immediate from (x), since $Cg(R)=R+R^-$. 

(xii) is proved by an induction similar to the one used in the proof of (iv).

(xiii) is immediate from (xii).

(xiv) Since
$Cg(R_1 \cup R_2 \cup \dots \cup R_n) =
R_1 + R _{1}^-+R_2 + R _{2}^-+\dots 
+ R _{n-1}  + R _{n-1}^- +R _{n} + R _{n}^-$, we get, by applying (xiii):
\begin{multline*} 
Cg(R_1 \cup R_2 \cup \dots \cup R_n) \subseteq
\\ 
\big(
F(R_1) + F(R _{1}^-) + F(R_2) + F(R _{2}^-) +\dots 
+F(R _{n} ) + F(R _{n}^-)
\big) \circ
\\
\overline{R_1 \cup R _{1}^-\cup R_2 \cup R _{2}^-\cup 
\dots \cup R _{n-1} \cup R _{n-1}^- \cup R _{n}\cup R _{n}^-} \circ
\\
\big(
G(R_1) + G(R _{1}^-) + G(R_2) + G(R _{2}^-) +\dots 
+G(R _{n} ) + G(R _{n}^-)
\big) 
\end{multline*} 
By (ix)  
we get 
$ R_i^- \subseteq F(R_i^-) \circ R_i \circ G(R_i^-)$, for all
$i=1,\dots,n$. Hence, 
$ R_1 \cup R _{1}^-\cup R_2 \cup R _{2}^-\cup 
\dots \cup R _{n-1} \cup R _{n-1}^- \cup R _{n}\cup R _{n}^-
\subseteq 
\big(
F(R _{1}^-) +  F(R _{2}^-) +\dots 
 + F(R _{n-1}^-) + F(R _{n}^-)
\big) \circ
(R_1 \cup R_2 \cup \dots \cup R _{n-1}  \cup R _{n})
\circ
\big(
G(R _{1}^-) +  G(R _{2}^-) +\dots 
 + G(R _{n-1}^-) + G(R _{n}^-)
\big)$.
Now (xiv) follows as in the proof of (x).
\end{proof} 

Condition (i) in Theorem \ref{wdta} has led us to 
the following result. If 
$F: \adma \to \adma$ let 
$F^{(2)}: \adma \to \adma$
be defined by $ F^{(2)}(R)=F(F(R))$.
Thus, in the next Theorem, 
$F^{(2)} \circ F^{(2)} \circ F^{(2)}$
is the operator $F'$ defined by
$F'(R)= F(F(R))\circ F(F(R)) \circ F(F(R))$,
and similarly for 
$G^{(2)} \circ G^{(2)} \circ G^{(2)}$. 

\begin{thm} \label{f2} 
If
$ {\alg A} $ has a term 
Mal$'$cev modulo $F$ and $G$
then $ {\alg A} $ has a term 
Mal$'$cev modulo $F^{(2)} \circ F^{(2)} \circ F^{(2)}$ and
$G^{(2)} \circ G^{(2)} \circ G^{(2)}$.
\end{thm}

We shall also be interested in the case when $F$ and $G$ are defined globally
on all algebras of some variety.

If $\v$ is a variety, let us say that $F$ is a
{\em  global operator on $\v$ for admissible and reflexive relations} 
if and only if to any algebra ${\alg A} \in \v$
$F$ assigns an operation $F_{\alg A}: \adma \to \adma $.
In case there is no danger of confusion, we shall omit the subscript ${\alg A} $.

We say that a global operator on $\v$
satisfies the {\em homomorphism property}
if and only if 
whenever
${\alg A}, {\alg  B} \in \v $,
$\phi: {\alg B} \to {\alg A} $
is a 
homomorphism, and 
$R \in \admb$ then
$\phi(F_{\alg B}(R)) \subseteq F_{\alg A}(\phi(R))$.
Here, $\phi(R)$ denotes the smallest compatible and
reflexive relation on 
${\alg A} $ which contains\
$\{(\phi (b), \phi(c))| bRc\} $.

As noticed in \cite[Remark 2.5]{lpcan2}, essentially all commutators
defined using matrices satisfy the homomorphism property. 

As usual, $F$ is said to be
{\em  monotone}
if and only if  $F(R) \subseteq F(S)$
whenever $R \subseteq S$.

\begin{thm} \label{wdtvar} 
Suppose that $\v$ is a variety,
 $F$, $G$ are 
global operators on $\v$ for admissible and reflexive relations, 
$F$, $G$ are monotone and satisfy the homomorphism property.
Then the following are equivalent:

(i)
$\v$ has a term which is 
Mal$'$cev modulo $F_{\alg A}$ and $G_{\alg A}$
for every algebra ${\alg A}$ in $\v$.

(ii) Every  $ {\alg A} \in \v$ has a term which is 
Mal$'$cev modulo $F_{\alg A}$ and $G_{\alg A}$.

(iii) The free algebra ${\alg X}$ in $\v$ generated by $2$ elements 
has a term which is 
Mal$'$cev modulo $F_{\alg X}$ and $G_{\alg X}$.

(iv) In the free algebra ${\alg X} $ in $\v$ generated by the two elements
$x$, $y $ there is a ternary term $t$ such that, if
$S$ is the smallest admissible and reflexive relation 
of ${\alg X}$
containing
$(x,y)$, then
\[
x F_{\alg X} (S) t(x,y,y) \qquad \text{ and } \qquad
t(x,x,y) G_{\alg X}(S) y
\]

(v) In every algebra ${\alg A}\in \v$ and 
for every 
relation $R \in \adma$, the following holds:
\[
R \circ R \subseteq F_{\alg A} (R) \circ R \circ G_{\alg A} (R)
\]

(vi) In the free algebra ${\alg X}$ in $\v$ generated by $3$ elements the following holds:
\[
R \circ R \subseteq F_{\alg X} (R) \circ R \circ G_{\alg X} (R)
\]
for every 
relation $R \in \admx$.

(vii) In the free algebra ${\alg X}$ in $\v$ generated by the three elements
$x$, $y$, $z$,  the following holds, where $S$ is the smallest admissible and reflexive relation of ${\alg X}$ containing
both $(x,y)$ and $(y,z)$:
\[
S \circ S \subseteq F_{\alg X} (S) \circ S \circ G_{\alg X} (S)
\]

(viii) In every algebra ${\alg A}\in \v$ and 
for every 
relation $R \in \adma$, the following holds:
\[
R \subseteq F_{\alg A} (R) \circ R^- \circ G_{\alg A} (R)
\]

(ix) In the free algebra ${\alg X}$ in $\v$ generated by $2$ elements the following holds:
\[
R \subseteq F_{\alg X} (R) \circ R^- \circ G_{\alg X} (R),
\]
for every 
relation $R \in \admx$.

(x) In the free algebra ${\alg X}$ in $\v$ generated by the two elements
$x$, $y$  the following holds, where $S$ is the smallest admissible and reflexive relation of ${\alg X}$ containing
$(x,y)$:
\[
S \subseteq F_{\alg X} (S) \circ S^- \circ G_{\alg X} (S),
\]
\end{thm} 

\begin{proof}
(i) $\Rightarrow $ (ii) $\Rightarrow $ (iii) $\Rightarrow $ (iv)
are trivial.

(iv) $\Rightarrow $ (i) 
We claim that the term $t$ given by (iv)
is 
Mal$'$cev modulo $F_{\alg A}$ and $G_{\alg A}$,
for every algebra ${\alg A}$ in $\v$.

Indeed, suppose that
${\alg A}$ in $\v$, 
$a, b \in {\alg A} $, 
$R \in \adma$,
and $aRb$.

Since ${\alg X} $ is
the free algebra  in $\v$ generated by 
$\{x,y\} $,
there is a homomorphism
$\phi: {\alg X} \to {\alg A} $
such that 
$\phi(x)=a$
and $\phi(y)=b$, and hence
$\phi (t(x,y,y))=t(a,b,b)$,
 $\phi (t(x,x,y))=t(a,a,b)$.

Since 
$x F_{\alg X} (S) t(x,y,y)$,
we have
$ \phi (x) \phi ( F_{\alg X} (S)) \phi ( t(x,y,y))$,
hence, by the homomorphism property,
$a F_{\alg A}(\phi (S)) t(a,b,b)$.

Since 
$S$ is the compatible and reflexive relation generated by $(x,y)$,
and $\phi$ is a homomorphism, 
then $\phi(S)$ is the compatible and reflexive relation generated by 
$(\phi(x),\phi(y))=(a,b)$,
and, since $a R b$, we have that 
$ \phi(S) \subseteq R$; thus, by the monotonicity of
$F_{\alg A}$, we get
$ F_{\alg A}(\phi(S)) \subseteq F_{\alg A}(R)$ and, eventually,
$a F_{\alg A}(R) t(a,b,b)$.

Exactly in the same way, we get
$t(a,a,b) G_{\alg A}(R) b$, thus
$t$ is 
Mal$'$cev modulo $F$ and $G$
for every algebra  in $\v$.

Having proved that
(iv) $\Rightarrow $ (i), 
we have that (i)-(iv) are all equivalent.

(i) $\Rightarrow$ (v) is from Theorem \ref{wdta}(vi).

(v) $\Rightarrow $ (vi) $\Rightarrow $ (vii) 
are trivial. 

(vii) $\Rightarrow $ (i). Since
$ xSySz $, we have 
$ (x,z) \in S \circ S$,
hence $ (x,z) \in  F_{\alg X} (S) \circ S \circ G_{\alg X} (S)$,
by assumption.

This means that ${\alg X}$ has terms $t_1(x,y,z)$ and  
$t_2(x,y,z)$ such that 
 $(x,
\spacex
 t_1(x,y,z)) \in F_{\alg X} (S)$,
$(t_1(x,y,z), t_2(x,y,z)) \in S$ and 
$( t_2(x,y,z),z) \in G_{\alg X} (S)$.

Notice that $ S= \{(u(x,y,z,x,y), u(x,y,z,y,z))| u \text{ a term of } 
{\alg X} \} $, since the right-hand relation is  reflexive,
admissible, and contains $(x,y)$ and $(y,z)$; moreover, every other reflexive admissible
relation containing $(x,y)$ and $(y,z)$ contains all pairs of the form
$(u(x,y,z,x,y), u(x,y,z,y,z))$.
Hence, $(t_1(x,y,z), t_2(x,y,z)) \in S$ means that 
there is a $5$-ary term $t'$ such that   
$t_1(x,y,z)=t'(x,y,z,x,y)$ and $t'(x,y,z,y,z)= t_2(x,y,z) $.

We claim that the ternary term 
$t(x,y,z)=t'(x,y,z,x,z)$ is Mal'cev modulo $F$ and $G$ throughout $\v$.
In order to prove it, first notice that   
$t_1(x,y,y)=t(x,y,y)$, and 
$t(y,y,z)=t_2(y,y,z)$.
Notice that, since ${\alg X}$ is a free algebra, 
the above identities hold throughout $\v$.

Suppose that
${\alg A} \in \v$, 
$a, b \in {\alg A} $, 
$R \in \adma$,
and $aRb$.
There is a homomorphism
$\phi: {\alg X} \to {\alg A} $
such that 
$\phi(x)=a$,
$\phi(y)=b$,
and $\phi(z)=b$, hence
$\phi (t_1(x,y,y))=t_1(a,b,b)$.

Since 
$x F_{\alg X} (S) t_1(x,y,z)$,
we have
$ \phi (x) \phi ( F_{\alg X} (S)) \phi ( t_1(x,y,z))$,
hence, by the homomorphism property,
$a F_{\alg A}(\phi (S)) t_1(a,b,b)$.

Since 
$S$ is the compatible and reflexive relation generated by $(x,y), (y,z)$,
then $\phi(S)$ is the compatible and reflexive relation generated by 
$(\phi(x),\phi(y))$, $(\phi(y),\phi(z))$, that is, generated by 
$(a,b), (b,b)$, hence generated simply by $(a,b)$.  
Since $a R b$, we have that 
$ \phi(S) \subseteq R$; thus, by the monotonicity of
$F_{\alg A}$, we get
$ F_{\alg A}(\phi(S)) \subseteq F_{\alg A}(R)$ and, eventually,
$a F_{\alg A}(R) t_1(a,b,b)= t(a,b,b)$.

Exactly in the same way,
by considering the homomorphism
$\psi: {\alg X} \to {\alg A} $
satisfying 
$\psi(x)=a$,
$\psi(y)=a$,
and $\psi(z)=b$, 
 we get
$t(a,a,b)=t_2(a,a,b) G_{\alg A}(R) b$.

Thus
$t$ is 
Mal$'$cev modulo $F$ and $G$
for every algebra  in $\v$, and
we have closed our second cycle of equivalencies.

(i) $\Rightarrow$ (viii) follows from Theorem \ref{wdta}(ix).

(viii) $\Rightarrow $ (ix) $\Rightarrow $ (x) 
are trivial. 

(x) $\Rightarrow $ (iv). Since $xSy$, then, by assumption,
$ (x,y) \in  F_{\alg X} (S) \circ S^- \circ G_{\alg X} (S)$.
This means that there are binary terms $t_1$ and $t_2$ such that 
 $ (x,t_1(x,y)) \in  F_{\alg X} (S) $,
$(t_1(x,y),t_2(x,y)) \in S^-$, and
$ (t_2(x,y),y) \in  G_{\alg X} (S) $.

That 
$(t_1(x,y),t_2(x,y)) \in S^-$ means that there is 
a ternary term $t$ such that 
$t_1(x,y)=t(x,y,y)$ and
$t(x,x,y)=t_2(x,y)$, since
$S= \{ (u(x,x,y),
\spacex
u(x,y,y))| u \text{ a term of } {\alg X}  \} $
 (cf., e.g., the proof of Theorem 1 (vi) $\Rightarrow$  (vii) in Part II,
or the proof of (vii) $\Rightarrow $ (i) here).

Thus,
$x  F_{\alg X} (S) t_1(x,y)= t(x,y,y)$,
$ t(x,x,y)=t_2(x,y) G_{\alg X} (S) y$,
that is, the hypotheses of (iv) are satisfied.
\end{proof}

\def\cprime{$'$} \def\cprime{$'$}
\providecommand{\bysame}{\leavevmode\hbox to3em{\hrulefill}\thinspace}
\providecommand{\MR}{\relax\ifhmode\unskip\space\fi MR }
\providecommand{\MRhref}[2]{%
  \href{http://www.ams.org/mathscinet-getitem?mr=#1}{#2}
}
\providecommand{\href}[2]{#2}

\end{document}